\documentclass[reqno, 12pt]{amsart}
\usepackage[top = 1in, bottom = 1in, left = 1in, right = 1in]{geometry}
\usepackage[hidelinks]{hyperref}
\usepackage{enumitem}

\DeclareMathOperator{\dive}{div}

\DeclareMathOperator{\trace}{trace}
\newtheorem{theorem}{Theorem}[section]

\begin{document}

\title[Interior Gradient Estimates for the Mean Curvature Equation]{Improved Interior Gradient Estimates for the Mean Curvature Equation under Nonlinear Assumptions}
\author[Xu]{Fanheng Xu}
\address{School of Mathematics (Zhuhai), Sun Yat-Sen University, 519082, Zhuhai, P. R. China}
\email{xufh7@mail.sysu.edu.cn}
\thanks{Xu was supported by the National Natural Science Foundation of China (No.12201654)}
\subjclass[2020]{Primary: 35B45, Secondary: 35J92, 35B50}
\keywords{Interior gradient estimate, Prescribed mean curvature equation, Maximum principle}

\begin{abstract}
In this paper, we investigate interior gradient estimates for solutions to the mean curvature equation
$$ \dive \left( \frac{\nabla u}{\sqrt{1 + |\nabla u|^2}} \right) = f(\nabla u)$$
under various nonlinear assumptions on the right-hand side. Under the weakened initial assumption $u\in C^1(B_R) \cap C^3(\{|\nabla u|>0\})$, we establish sharp gradient bounds that depend on the oscillation of the solution. These estimates are applicable to a wide class of nonlinear terms, including the specific forms arising from the elliptic regularization of the inverse mean curvature flow ($f=\varepsilon\sqrt{1+|\nabla u|^2}$ ), minimal surface equation ($f=0$) and several polynomial and logarithmic growth regimes. As applications, the gradient bounds imply uniform ellipticity of the equation away from the critical set,which allows one to apply classical elliptic regularity theory and obtain higher regularity of solutions in the noncritical region. Moreover, when the solution grows at most linearly, all cases of our results can be applied in Moser's theory to establish the affine linear rigidity of global solutions. This directly leads to the Liouville-type theorems for global solutions without requiring additional proofs.
\end{abstract}

\maketitle

\section{Introduction}

Gradient estimates play a crucial role in the theory of elliptic equations, offering critical insights into the behavior of solutions, including existence, regularity, and uniqueness results.
Let $u_{ij} = u_{x_i x_j}$, and adopt the summation convention.
For the quasilinear equation of the following form
\begin{align}\label{Q}
\mathcal{Q}u = a^{ij} (x, u, \nabla u) u_{ij} + b(x, u, \nabla u) = 0,
\end{align}
if $a^{ij}$ and $b$ satisfy certain structural conditions such that operator $\mathcal{Q}$ is uniformly elliptic, then the gradient of any solution $u$ to (\ref{Q}) in a ball $B_R(0)$ of $\mathbb{R}^n$ usually satisfies
\begin{align}\label{std uniform elliptic grad est}
|\nabla u(0)| \leq C\left( 1 + \frac{1}{R^{\frac{1}{\theta}}} \right).
\end{align}
For specific structural conditions, one may refer to \cite{Gilbarg Trudinger Berlin 1998} and \cite{Ladyzhenskaya Uraltseva book 1968}.
For the non-uniformly elliptic case, such as the prescribed mean curvature equation
\begin{align}\label{prescribed mean curvature equation}
\dive\left( \frac{\nabla u}{\sqrt{1 + |\nabla u|^2}} \right) = H(x).
\end{align}
The gradient of any solution $u$ to (\ref{prescribed mean curvature equation}) in $B_R(0)$ satisfies
\begin{align}\label{std mean curvature eq grad est}
|\nabla u(0)| \leq \exp\left( C \left( 1 + \frac{\sup_{B_R(0)}(u-u(0))}{R} \right) \right).
\end{align}

The earliest gradient estimate (\ref{std mean curvature eq grad est}) was established by Finn \cite{Finn ARMA 1963} for dimension $n=2$ and later generalized by Bombieri, De Giorgi and Miranda \cite{Bombieri Giorgi Miranda ARMA 1969} for arbitrary $n$. Their proofs relied on properties of minimal surfaces, such as the isoperimetric inequality and the resulting Sobolev inequality.
Since then, many scholars have proposed various methods for proving gradient estimates similar to (\ref{std mean curvature eq grad est});
see
\cite{Bombieri Giusti Invent 1972}
\cite{DeSilva Jerison CAPM 2011}
\cite{Ladyzhenskaya Uraltseva CPAM 1970}
\cite{Wang MathZ 1998}.
In particular, Trudinger \cite{Trudinger PNAS 1972} employed classical potential theory, which can be extended to derive gradient bounds for arbitrary $C^2$ functions in terms of the mean curvature and the gradient of the mean curvature of the graph.
Korevaar \cite{Korevaar PSPM 1983} derived the gradient estimate using the technique of normal variations. In his approach, the mean curvature of the surface varied in the normal direction is expressed in terms of the curvature of the original surface. By employing an appropriate test function, the priori estimate is obtained.
The Dirichlet problem for the prescribed mean curvature equation was studied by
Serrin \cite{Jenkins Serrin 1968 JRAM, Serrin PTRSLA 1969}
and
L. Simon \cite{Simon Indiana 1976}.
For the Neumann boundary value problem, please refer to Ma and Xu \cite{Ma Xu Adv 2016} and the references therein.
A more detailed history could be found in Gilbarg and Trudinger \cite{Gilbarg Trudinger Berlin 1998}.

By closely examining these two gradient estimates (\ref{std uniform elliptic grad est}) and (\ref{std mean curvature eq grad est}), we observe that the upper bound of the gradient estimate always remains larger than a positive constant as $R \to \infty$.
Nevertheless, we can find examples where the global solution has vanishing gradients.

For uniformly elliptic equations, an example is the classic Liouville theorem, which states that all bounded entire harmonic functions, i.e., solutions to the Laplace equation $\Delta u=0$ must be constant, implying that its gradient vanishes;
For the non-uniformly elliptic case, an example is the minimal surface equation
\begin{align}\label{minimal surface equation}
\dive \left( \frac{\nabla u}{\sqrt{1 + |\nabla u|^2}} \right) = 0.
\end{align}
Bernstein \cite{Bernstein MathZ 1927} first proved that, in two dimensions $n = 2$, any global solution to the minimal surface equation (\ref{minimal surface equation}) in $\mathbb{R}^n$ must be an affine linear function.
Later, Fleming \cite{Fleming RCMP 1962}, De Giorgi \cite{DeGiorgi ASNSP 1965}, Almgren \cite{Almgren Annals 1966}, and J. Simons \cite{Simons Annals 1968} extended Bernstein's result to dimensions $n \leq 7$.
For dimensions $n \geq 8$, Bernstein's theorem no longer holds, since Bombieri et al. \cite{Bombieri Giorgi Giusti Invent 1969} found a counterexample.
However, Moser's series of results in \cite{Moser CPAM 1961} indicate that if we add the condition that the gradient of the global solution is bounded, then global solution remains affine linear functions. This condition is satisfied by the gradient estimate (\ref{std mean curvature eq grad est}) when the solution $u$ grows at most linearly on one side, meaning
\begin{align*}
u(x) \leq C(1 + |x|), \quad \forall x\in \mathbb{R}^n.
\end{align*}
To summarize, if the growth of the global solution of (\ref{minimal surface equation}) is constrained linearly on one side, then the affine linear rigidity of the global solution can be established for all dimensions (More rigidity discussion could be found in \cite{Ecker Huisken JDG 1990} and \cite{Simon JDG 1989}).
Furthermore, if we additionally assume that the global solution is also bounded on one side, we can show the global solution to the minimal surface equation (\ref{minimal surface equation}) must be constant.
These examples illustrate that for a wide class of quasilinear elliptic equations---both uniformly and non-uniformly elliptic---the gradient estimates for bounded solutions may tend to zero as the domain radius increases.

This leads to an interesting question: whether, under certain conditions, some constants in the estimates (\ref{std uniform elliptic grad est}) and (\ref{std mean curvature eq grad est}) can be removed to improve the gradient bounds,
so that the gradient upper bound tends to $0$ as $R\rightarrow\infty$
(under the assumption that the solution is bounded or satisfies some sublinear growth conditions)?

For the uniformly elliptic case, it is well-known that solutions to the Laplace equation satisfy gradient estimate
\begin{align*}
|\nabla u(0)| \leq n \frac{\sup_{B_R(0)} |u|}{R}
\end{align*}
This provides an affirmative answer to this question.
However, for the non-uniformly elliptic case, even for the minimal surface equation (\ref{minimal surface equation}), no analogous results have been established.

We aim to study the mean curvature equation of the following form
\begin{align}\label{eq}
\dive \left( \frac{\nabla u}{\sqrt{1 + |\nabla u|^2}} \right) = f(\nabla u)
\quad \text{in} \quad \mathbb{R}^n.
\end{align}
where $f=f(p)$ denotes a $C^1$ function of $p\in\mathbb{R}^n$.
We impose different assumptions on the nonlinear term $f$ to derive the desired gradient estimates.
Let
\begin{align*}
\begin{cases}
\nabla_p f = (f_{p_1},f_{p_2},\ \cdots, f_{p_n}),\\
|\nabla_p f| = \left(\sum_{j=1}^n |f_{p_j}|^2\right)^{1/2},\\
L := \sup\limits_{x, y \in B_R(0)} |u(x)-u(y)|.
\end{cases}
\end{align*}
We first introduce the following conditions depending on the growth of $f$:
\begin{enumerate}
\item[(A1)] For $\theta > 0$, there exist constants $m_1, m_2 > 0$ such that
\[
f(p)^2 - m_1 |p|^2 |\nabla_p f(p)|^2 \geq m_2 |p|^{2\theta}.
\]

\item[(A2)] For $\theta > 1$, there exist constants $m_1, m_2, m_3 > 0$ such that
\[
\begin{cases}
f(p)^2 - m_1 (1 + |p|^2) \log(1 + |p|^2) |\nabla_p f(p)|^2 \geq m_2\log^{2\theta}(1 + |p|^2), \\
|f(p)| \leq m_3 \log^{\theta}(1 + |p|^2).
\end{cases}
\]

\item[(A3)] For $0 < \theta \leq 1$, there exists $m_1 > 0$ such that
\[
|f(p)| = m_1 \log^{\theta}(1 + |p|^2).
\]

\item[(A4)] There exist constants $m_1, m_2 > 0$ such that
\[
\begin{cases}
f(p)^2 - m_1 (1 + |p|^2)^2 \log(1 + |p|^2) |\nabla_p f(p)|^2 \geq 0, \\
|f(p)| \leq m_2.
\end{cases}
\]
\end{enumerate}

Our main result is the following:
\begin{theorem}\label{main theorem}
Let $u\in C^1(B_R(0)) \cap C^3(B_R(0)\cap\{|\nabla u|>0 \})$ be a solution of (\ref{eq}).

\begin{itemize}
\item[(a)] Suppose $\theta > 0$ and $f$ satisfies (A1). Then $u$ satisfies
\begin{align}\label{thm est 1}
|\nabla u(0)| \leq \frac{C}{R^{\frac{1}{\theta}}},
\end{align}
where $C$ depends on $n$, $m_1$, $m_2$ and $\theta$.

\item[(b)] Suppose $\theta > 1$ and $f$ satisfies (A2). Then $u$ satisfies
\[
|\nabla u(0)|^2
\leq C \max\left\{\exp\left(\frac{1}{R^{\frac{1}{\theta-1}}} \right)-1,
\ \frac{1}{R^{\frac{2}{2\theta-1}}}\right\},
\]
where $C$ depends on $n$, $m_1$, $m_2$, $m_3$ and $\theta$.

\item[(c)] Suppose $0 < \theta \leq 1$ and $f$ satisfies (A3). Then for any $0<\eta<1$, $u$ satisfies
\[
|\nabla u(0)|^2 \leq \max\left\{\exp\left(\frac{C_1(L+1)^{\frac{1}{2\theta}}}{R^{\frac{1}{\theta}}}\right)-1,\
\exp\left(\frac{C_2(L+1)^{\frac{1}{1-\eta}}}{R^{\frac{2}{1-\eta}}}\right)-1\right\},
\]
where $C_1$ depends on $n$, $m_1$, $\theta$, and $C_2$ depends on $n$, $m_1$, $\eta$, $\theta$.

\item[(d)] Suppose $f$ satisfies (A4), then $u$ satisfies
\[
|\nabla u(0)|^2 \leq \max\left\{\exp\left(\frac{C(L+1)^2}{R^2}\right)-1,\
\exp\left(\frac{C(L+1)^2}{R^{\frac{2}{3}}}\right)-1\right\},
\]
where $C$ depends on $n$, $m_1$, $m_2$.

\item[(e)] Suppose $f \equiv 0$, then $u$ satisfies
\[
|\nabla u(0)|^2 \leq \max\left\{\exp\left(\frac{CL^2}{R^2}\right)-1,\
\exp\left(\frac{CL}{R}\right)-1\right\},
\]
where $C$ depends only on $n$.
\end{itemize}
\end{theorem}

Here are some comments on Theorem \ref{main theorem}.

\begin{itemize}[leftmargin=*]

\item
Theorem \ref{main theorem} can be extended to the case $f = f(x, v, p) \in C^1( \mathbb{R}^n \times \mathbb{R} \times \mathbb{R}^n )$, where $v=u$ and $p = \nabla u$.
Some special cases include:
\begin{itemize}
\item
When $f=f(x,p)$ and $|f_x| \leq C$,
under the conditions in Theorem \ref{main theorem} (a), we can obtain appropriate gradient estimates for $\theta>1/2$.
\item
When $f = f(v, p)$ and $f_v\geq 0$, all cases in Theorem \ref{main theorem} can be applied to such $f$.
\end{itemize}
It can be observed that our theorem indicates that the gradient upper bound is influenced significantly more by the order of $f$ with respect to $p$ than by its dependence on $x$ or $v$.
That is why we assume that nonlinear term $f$ depends only on the gradient of $u$ in the equation (\ref{eq}).
Another reason is that previous works have primarily focused on
$f=0$ (please refer to \cite{Bombieri Giorgi Miranda ARMA 1969} \cite{Finn ARMA 1963}) or $f=f(x,v)$ (please refer to \cite{Korevaar PSPM 1983}\cite{Ma Xu Adv 2016}\cite{Wang MathZ 1998}), whereas this paper addresses the case of $f=f(p)$ filling a gap left in previous literature.

\item
The mean curvature operator in (\ref{eq}) can also be extended to more general quasilinear operators of the form similar to (\ref{Q}) using the approach of Serrin \cite{Serrin PTRSLA 1969}.
In equation (\ref{Q}), we may consider
\begin{align*}
a^{ij}(x, v, p) = a^{ij}_*(p) + \frac{1}{2} \left( p_ic_j(x, v, p) + c_i(x, v, p)p_j \right).
\end{align*}
Here $a^{ij}_*(p)\in C^1(\mathbb{R}^n)$ and $c_i\in C^1(\Omega\times\mathbb{R}\times\mathbb{R}^n)$, where $i$, $j = 1, \cdots, n$.
For the equation (\ref{eq}) under consideration, see also its reformulated version (\ref{eq rewrite}), we have
\begin{align}\label{f714dbbcc}
a^{ij}_* = \delta^{ij}, \qquad
c_i = -\frac{p_i}{1 + |p|^2}.
\end{align}
By simply imposing some structure conditions on $a^{ij}_*$ and $c_i$ similarly to (\ref{f714dbbcc}), our proof can also be extended to general quasilinear equations of the form (\ref{Q}). Specific implementations of this approach can also be found in Serrin \cite{Serrin JFA 1970} and in the book of Ladyzhenskaya and Ural'tseva \cite{Ladyzhenskaya Uraltseva book 1968}, and Gilbarg and Trudinger \cite{Gilbarg Trudinger Berlin 1998}.
For further approaches to quasilinear equations, please refer to \cite{Caffarelli Garofalo Segala CPAM 1994} \cite{Cozzi Farina Valdinoci CMP 2014} \cite{Farina Valdinoci CVPDE 2014}.
Additionally, for approaches to fully nonlinear equations, please refer to \cite{Chou Wang CPAM 2001}.

\item
Regarding the minimal surface equation ($f=0$),
we prove in Theorem \ref{main theorem}(e) a gradient estimate stronger than the one by Bombieri-De Giorgi-Miranda,
which directly implies that minimal graphs over $\mathbb{R}^n$ growing less than linearly on one side are constant.
Furthermore, our results are applicable to a wide class of interesting nonlinear terms $f(p)$. Notable examples include:
\begin{itemize}
\item
$f(p) = |p|^\theta$, which satisfies (A1) with $\theta>0$.
In particular, the case $\theta = 1$ is of special geometric interest as it covers the equation
\begin{align*}
\text{div}\left(\frac{\nabla u}{\sqrt{1+|\nabla u|^2}}\right) = \varepsilon \sqrt{1+|\nabla u|^2};
\end{align*}
arising from the elliptic regularization of the inverse mean curvature flow (set $u:=u_\varepsilon /\varepsilon $ in \cite{Huisken Ilmanen JDG 2001}) studied by Huisken and Ilmanen.
\item
$f(p) = \log^\theta(1 + |p|^2)$,
which satisfies (A2) with $\theta>1$ and (A3) with $0<\theta\leq 1$,
representing the case of logarithmic growth;
\item
$f(p) = \frac{|p|}{\sqrt{1 + |p|^2}}$, which satisfies (A4), representing the case of no growth.
\end{itemize}

\item
The most recent result directly related to our work is found in \cite{Wang MathZ 1998}.
In the following, we will compare our methods and results with those in \cite{Wang MathZ 1998}.
Both approaches adopt the idea of controlling the gradient term through structural analysis of the operator at a single point.
However, by establishing a lower bound estimate on the Hessian of the solution, we successfully achieve control over both the nonlinear terms $f$ and $f^2$,
rather than only controlling $f$ as in \cite{Wang MathZ 1998}.
This is crucial when $f$ involves higher-order terms of $\nabla u$,
such as $f=|\nabla u|^\theta$.
This advancement allows us to derive estimates of the form (\ref{std uniform elliptic grad est}), as demonstrated in (\ref{thm est 1}) and its proof.
We also improve the auxiliary function.
After studying the auxiliary function
$\log(|\nabla u|)$ used in \cite{Wang MathZ 1998} and
$\log(|\nabla u|^2)$ used in \cite{Ma Xu Adv 2016},
we find a more effective auxiliary function $\log(1+|\nabla u|^2)$ for the mean curvature operator.
By a careful calculation on this auxiliary function, an upper bound for the gradient is obtained.

\item
As an application of the interior gradient estimates in Theorem \ref{main theorem},
we note that the local boundedness of $|\nabla u|$ plays a fundamental role in the regularity theory.
In the settings of Theorem \ref{main theorem}, provided that $u$ satisfies the prescribed growth conditions (specifically in cases (c), (d), and (e)),
the gradient satisfies $|\nabla u| \leq C$ on any compact subdomain of $B_R(0)$.
Consequently, the coefficient matrix
$$a^{ij}(p) = \delta^{ij} - \frac{p_i p_j}{1 + |p|^2}$$
of mean curvature operator becomes uniformly elliptic, meaning there exist positive constants $\lambda, \Lambda$ such that
$$\lambda |\xi|^2 \leq a^{ij}(\nabla u) \xi_i \xi_j \leq \Lambda |\xi|^2$$
for all $\xi \in \mathbb{R}^n$.
This uniform ellipticity allows us to invoke classical elliptic regularity theory,
such as Schauder estimates or $W^{2,p}$ theory (see \cite{Gilbarg Trudinger Berlin 1998}),
to bootstrap the regularity of solutions from the initial $C^1$ assumption to $C^\infty$ (assuming $f$ is smooth) in the region where $|\nabla u| > 0$.
Given the local boundedness of $|\nabla u|$ and then the non-degenerate nature of the mean curvature operator, standard regularity theory for quasilinear elliptic equations further ensures that this smoothness extends across the entire ball $B_R(0)$.

\item
As a further application of Theorem \ref{main theorem}, we derive Liouville-type theorems for global solutions to equation \eqref{eq} in $\mathbb{R}^n$.
Under the respective assumptions (A1)--(A4), any entire solution must be constant in the following cases:
\begin{enumerate}
\item[(E1)] $f(p) = |p|^\theta$ with $\theta > 0$;
\item[(E2)] $f(p) = \log^\theta(1 + |p|^2)$ with $\theta > 1$;
\item[(E3)] $f(p) = \log^\theta(1 + |p|^2)$ with $0 < \theta \leq 1$, provided $u(x) = o(|x|^2)$ as $|x| \to \infty$;
\item[(E4)] $f(p) = \frac{|p|}{\sqrt{1 + |p|^2}}$, provided $u(x) = o(|x|^{1/3})$ as $|x| \to \infty$.
\end{enumerate}
In comparison, similar rigidity results are obtained in
\cite{Bianchini Colombo Magliaro Mari Pucci Rigoli 2021 MathEng} and
\cite{Bianchini Mari Pucci Rigoli 2021 Springer} using different methods, including maximum principles at infinity and a Keller-Osserman theory for quasilinear inequalities on complete manifolds. For instance, \cite[Theorem 2.21]{Bianchini Colombo Magliaro Mari Pucci Rigoli 2021 MathEng} establishes constancy for (E4) without the $o(|x|^{1/3})$ growth condition, and similarly for (E1) and (E3) with $\theta \in (0,1/2]$ in $\mathbb{R}^n$ without growth assumptions.
However, these methods do not cover cases like (E1) and (E2) for $\theta > 1/2$,
where our Theorem \ref{main theorem} provides a unified framework that not only covers the rigidity but also yields a precise quantitative decay of the gradient as $R \to \infty$. This quantitative nature is particularly useful when dealing with more general nonlinear terms $f(p)$ where the monotonic methods of \cite{Bianchini Colombo Magliaro Mari Pucci Rigoli 2021 MathEng} and
\cite{Bianchini Mari Pucci Rigoli 2021 Springer} might not be directly applicable.

\end{itemize}

The method we use in the proof is a combination of the maximum principle and the Bernstein technique.
The maximum principle is a effective tool for deriving gradient estimates. By establishing the upper (or lower) harmonicity of the gradient function and comparing its boundary values with its internal maximum, gradient bounds can be obtained. This approach is particularly useful when the boundary conditions of the solution are known.
Bernstein's method, originally introduced by Bernstein himself in \cite{Bernstein MathAnn 1906} and \cite{Bernstein MathAnn 1910}, is a powerful non-variational tool with applications across a wide range of equations and contexts.
In Bidaut-V\'eron et al. \cite{Bidaut-Veron Gacia-Huidobro Veron Duke 2019}, the authors
introduced highly effective auxiliary functions, and successfully applied Bernstein's method to quasilinear differential equations in both pointwise and integral forms, obtaining corresponding gradient estimates.
Recently, Cabr\'e, Dipierro and Valdinoci \cite{Cabre Dipierro Valdinoci ARMA 2022} extended Bernstein's method to the context of integro-differential operators. They derived first-order and one-sided second-order derivative estimates for solutions to fractional equations.

The aim of the following section is to prove Theorem \ref{main theorem}. The proof is divided into six steps. In the first step, we prepare some preliminary work for the proof.
In Step 2 through 6, we complete the proof of Theorem \ref{main theorem} (a) through (e), respectively.

\section{Interior Gradient Estimates}
Here and in the sequel, we denote $z = |\nabla u|^2$ and use $\delta^{ij}$ to denote the Kronecker delta symbol.
Unless otherwise specified, we use subscripts to denote partial derivatives
\begin{align*}
(\cdot)_i := \partial_{x_i}(\cdot), \quad (\cdot)_{ij} := \partial_{x_ix_j} (\cdot), \quad \text{etc}.
\end{align*}
We employ the summation convention, meaning that if the same index appears twice within a single term, it is implicitly summed over and thus eliminated.
The letters $C$, $C_{n,\alpha}$, etc., denote positive constants whose values are unimportant and may vary at different occurrences.

\subsection*{Step 1.}
The equation (\ref{eq}) can be rewritten as
\begin{align}\label{eq rewrite}
a^{ij}u_{ij} := \left( \delta^{ij} -\frac{u_iu_j}{1 + z} \right)u_{ij} = g := (1 + z)^\frac{1}{2}f.
\end{align}
At point where $|\nabla u|>0$,
above equation can be differentiated
with respect to $x_k$, yielding
\begin{align*}
\left( -\frac{u_{ik}u_j + u_iu_{jk}}{1 + z} + \frac{u_iu_jz_k}{(1 + z)^2} \right)u_{ij} + a^{ij} u_{ijk} = \frac{\mathrm{d} g}{\mathrm{d} x_k}, \quad k = 1, 2, \cdots, n.
\end{align*}
Multiplying both sides by $u_k$ and summing over the index $k$, we derive
\begin{align}\label{3to2}
a^{ij}u_{ijk}u_k = \frac{\mathrm{d}g}{\mathrm{d}x_k} u_k + \frac{|\nabla z|^2}{2(1 + z)} - \frac{(\nabla z \nabla u)^2}{2(1 + z)^2}.
\end{align}
Here, we have used that $u_k u_{jk} = z_j/2$.

Let $h = h(u)$ and $F = F(z)$ be some functions to be determined, satisfying
\begin{align*}
h>0, \quad F>0, \quad F'>0.
\end{align*}
Let $\alpha\geq 1$ be a sufficiently large constant to be determined.
Define $\varphi$ as a cut-off function of the following form:
\begin{align}\label{def varphi}
\varphi(x) := \left( 1-\frac{|x|^2}{R^2} \right)^\alpha.
\end{align}
We consider the maximum point of $hF\varphi$ in $B_R(0)$, denoted by $x_0$.
We may assume that $\nabla u \neq 0$ at point $x_0$ otherwise gradient estimate is trivial.
Unless otherwise specified, all calculations will be performed at the point $x_0$.

Then at point $x_0$, there holds
\begin{align}\label{87e935eb9}
0 & \geq a^{ij} \left( \log(hF\varphi) \right)_{ij} \notag\\
& = a^{ij} \left( \left( \frac{h''}{h} - \frac{h'^2}{h^2} \right) u_iu_j + \frac{h'}{h}u_{ij} + \left( \frac{F''}{F}-\frac{F'^2}{F^2} \right) z_iz_j + \frac{F'}{F} z_{ij} + \frac{\varphi_{ij}}{\varphi} - \frac{\varphi_i\varphi_j}{\varphi^2} \right).
\end{align}

We now compute each term on the right-hand side of the above inequality.
For the terms involving $u_i$ and $u_{ij}$, we have
\begin{align}\label{est I1}
a^{ij} \left( \left( \frac{h''}{h} - \frac{h'^2}{h^2} \right) u_iu_j + \frac{h'}{h}u_{ij} \right) = \left( \frac{h''}{h} - \frac{h'^2}{h^2} \right) \frac{z}{1 + z} + \frac{h'}{h} g.
\end{align}
For the terms involving $z_i$, we have
\begin{align}\label{est I21 I22}
a^{ij}\left( \frac{F''}{F}-\frac{F'^2}{F^2} \right) z_iz_j = \left( \frac{F''}{F}-\frac{F'^2}{F^2} \right)\left( |\nabla z|^2 - \frac{(\nabla z\nabla u)^2}{1 + z} \right).
\end{align}
Using $u_k u_{jk} = z_j/2$ again, we get
\begin{align}\label{c6db997b9}
a^{ij} u_{ik} u_{jk} = \|\nabla^2 u\|^2 - \frac{|\nabla z|^2}{4(1 + z)},
\end{align}
where $\nabla^2 u$ is the Hessian matrix of $u$, and $\|\nabla^2 u\|^2 := \sum_{i, j} u_{ij}^2$ is the squared Hilbert-Schmidt norm of $\nabla^2 u$.
Then for the term involving $z_{ij}$, using (\ref{c6db997b9}) and (\ref{3to2}), we obtain
\begin{align}\label{8295406c3}
\frac{F'}{F} a^{ij} z_{ij} & = \frac{F'}{F}a^{ij} \left( 2u_{jk}u_{ik} + 2u_ku_{ijk} \right) \notag\\
& = \frac{F'}{F} \left( 2\|\nabla^2 u\|^2 + \frac{|\nabla z|^2}{2(1 + z)}
+ 2 \frac{\mathrm{d} g}{\mathrm{d} x_k} u_k - \frac{(\nabla z \nabla u)^2}{(1 + z)^2}
\right).
\end{align}
Let $\bf A$ be the $n\times n$ matrix with $(i, j)^{th}$ entry $a^{ij}$, then we have
\begin{align*}
\trace({\bf A} \nabla^2 u) = a^{ij} u_{ij} = g
\end{align*}
and
\begin{align*}
\|{\bf A} \nabla^2 u\|^2 & = \sum_{i, j} \left( a^{ik} u_{jk} \right)^2 \notag\\
& = \sum_{i, j} \left( u_{ij} - \frac{u_i z_j}{2(1 + z)} \right)^2 \notag\\
& = \sum_{i, j} \left( u_{ij}^2 - \frac{u_{ij} u_i z_j}{1 + z} + \frac{u_i^2 z_j^2}{4(1 + z)^2} \right)\notag\\
& = \|\nabla^2 u\|^2 - \frac{2 + z}{4(1 + z)^2}|\nabla z|^2.
\end{align*}
Combining these with the Cauchy inequality
\begin{align*}
\frac{1}{n} \left( \trace({\bf A} \nabla^2 u) \right)^2 \leq \|{\bf A} \nabla^2 u\|^2,
\end{align*}
we derive a lower bound for $\nabla^2 u$, given by
\begin{align*}
\|\nabla^2 u\|^2\geq \frac{1}{n} g^2 + \frac{2 + z}{4(1 + z)^2}|\nabla z|^2.
\end{align*}
Substituting the above into (\ref{8295406c3}), we obtain
\begin{align}\label{est I23}
\frac{F'}{F} a^{ij} z_{ij} \geq \frac{F'}{F} \left( 2\frac{\mathrm{d} g}{\mathrm{d} x_k}u_k + \frac{2}{n}g^2 + \frac{3 + 2z}{2(1 + z)^2}|\nabla z|^2 - \frac{(\nabla z\nabla u)^2}{(1 + z)^2} \right).
\end{align}

Recalling the definition of $\varphi$ in (\ref{def varphi}), we have
\begin{align}\label{varphi derivative}
\begin{cases}
\varphi_i = -2\alpha\varphi^{\frac{\alpha-1}{\alpha}}\frac{x_i}{R^2}; \\
\varphi_{ij} = 4\alpha(\alpha-1) \varphi^{\frac{\alpha-2}{\alpha}} \frac{x_ix_j}{R^4} - 2\alpha\varphi^\frac{\alpha-1}{\alpha} \frac{\delta_{ij}}{R^2}. \\
\end{cases}
\end{align}
This implies that, for the term involving $\varphi_i$ and $\varphi_{ij}$,
note that $|\varphi|<1$,
we have the estimate
\begin{align}\label{est I3}
a^{ij}\left( \frac{\varphi_{ij}}{\varphi} - \frac{\varphi_i\varphi_j}{\varphi^2} \right)
\geq - \frac{C_{n,\alpha}}{R^2 \varphi^\frac{2}{\alpha}},
\end{align}
where $C_{n,\alpha}$ is a positive constant depending on $n$ and $\alpha$.
Combining (\ref{est I1}), (\ref{est I21 I22}), (\ref{est I23}), and (\ref{est I3}) with (\ref{87e935eb9}), we conclude that
\begin{align}\label{ad6b8fa97}
\frac{F'}{F} \left( 2\frac{\mathrm{d} g}{\mathrm{d} x_k}u_k + \frac{2}{n}g^2 \right)
+ \frac{h'}{h} g
+ \left( \frac{h''}{h} - \frac{h'^2}{h^2} \right) \frac{z}{1 + z} \notag\\
+ G |\nabla z|^2 + H \frac{(\nabla z \nabla u)^2}{1 + z}
\leq \frac{C_{n,\alpha}}{R^2\varphi^\frac{2}{\alpha}},
\end{align}
where
\begin{align}\label{def G H}
\begin{cases}
\displaystyle
G := \frac{F''}{F} - \frac{F'^2}{F^2} + \frac{3 + 2z}{2(1 + z)^2}\frac{F'}{F}; \\
\displaystyle
H := -\frac{F''}{F} + \frac{F'^2}{F^2} - \frac{1}{1 + z} \frac{F'}{F}.
\end{cases}
\end{align}

Recalling that $x_0$ is the maximum point, we also have at $x_0$
\begin{align*}
\nabla \left( hF\varphi \right) = 0.
\end{align*}
This is equivalent to
\begin{align}\label{s1 est nabla z}
\nabla z = -\frac{F}{F'} \left( \frac{\nabla \varphi}{\varphi} + \frac{h'}{h} \nabla u \right).
\end{align}
Then, using (\ref{s1 est nabla z}), we have
\begin{align}\label{fbfcc84c5}
&G |\nabla z|^2 + H \frac{(\nabla z \nabla u)^2}{1 + z} \notag\\
& = \left( \frac{h'F}{hF'} \right)^2 \left( G z + H \frac{z^2}{1 + z} \right)
+ \frac{2h'F^2}{hF'^2\varphi} \left(G + H \frac{z}{1 + z} \right) \nabla u\nabla \varphi \notag\\
&\qquad + \left( \frac{F}{F'\varphi} \right)^2\left( G |\nabla\varphi|^2 + H \frac{(\nabla u \nabla\varphi)^2}{1 + z} \right).
\end{align}
Recalling that $g = (1 + z)^{\frac{1}{2}} f$ and applying (\ref{s1 est nabla z}) once again, we have
\begin{align*}
\frac{\mathrm{d} g}{\mathrm{d} x_k} & = (1 + z)^{\frac{1}{2}} \frac{\mathrm{d} f}{\mathrm{d} x_k} + \frac{z_k f}{2(1 + z)^{\frac{1}{2}}} \notag\\
& = \sum_j (1 + z)^{\frac{1}{2}} f_{p_j} u_{jk}
- \frac{F f \varphi_k}{2F'(1 + z)^{\frac{1}{2}}\varphi}
- \frac{F h' f u_k}{2F'h(1 + z)^{\frac{1}{2}}}.
\end{align*}
Then we have
\begin{align}\label{926b5bd10}
& 2\frac{F'}{F} \frac{\mathrm{d} g}{\mathrm{d} x_k} u_k + \frac{h'}{h} g \notag\\
& = \sum_j \frac{F'}{F}(1 + z)^{\frac{1}{2}} f_{p_j} z_j
- \frac{f \nabla u \nabla \varphi}{(1 + z)^{\frac{1}{2}} \varphi}
+\frac{h'f}{h(1 + z)^{\frac{1}{2}}}.
\end{align}
Substituting (\ref{fbfcc84c5}) and (\ref{926b5bd10}) into (\ref{ad6b8fa97}), we derive the following inequality
\begin{align}\label{step1 result}
I_1 + I_2 + I_3 + I_4 + I_5 \leq I_6 + I_7 + I_8 + I_9,
\end{align}
where
\begin{align*}
&I_1 := \sum_j \frac{F'}{F}(1 + z)^{\frac{1}{2}} f_{p_j} z_j, \\
&I_2 := \frac{2}{n}\frac{F'}{F}(1 + z)f^2, \\
&I_3 := \frac{h'f}{h(1 + z)^{\frac{1}{2}}}, \\
&I_4 := \left( \frac{h''}{h} - \frac{h'^2}{h^2} \right) \frac{z}{1 + z}, \\
&I_5 := \left( \frac{h'F}{hF'} \right)^2 \left( G z + H \frac{z^2}{1 + z} \right), \\
&I_6 := -\frac{2h'F^2}{hF'^2\varphi} \left(G + H \frac{z}{1 + z} \right) \nabla u\nabla \varphi, \\
&I_7 := \frac{f \nabla u \nabla \varphi}{(1 + z)^{\frac{1}{2}} \varphi}, \\
&I_8 := \frac{F^2}{F'^2\varphi^2}\left( -G |\nabla\varphi|^2 - H\frac{(\nabla u\nabla \varphi)^2}{1 + z} \right), \\
&I_9 := \frac{C_{n,\alpha}}{R^2 \varphi^{\frac{2}{\alpha}}}.
\end{align*}

In the subsequent steps, we will choose appropriate $h$ and $F$, and impose different conditions on $f$ to obtain the desired gradient estimates.

\subsection*{Step 2.}
In this step, we continue from Step 1 to prove Theorem \ref{main theorem} (a).
Although there are alternative approaches to handle the terms involving $h$ in (\ref{step1 result}) (as discussed in the subsequent steps), the simplest way in this step is to set $h=1$. This eliminates all terms containing $h'$ and $h''$ in (\ref{step1 result}), leading to
\begin{align}\label{6e8b74b95}
I_1 + I_2 \leq I_7 + I_8 + I_9.
\end{align}
By further setting $F = z$, (\ref{s1 est nabla z}) simplifies to
\begin{align*}
\nabla z = - z \frac{\nabla \varphi}{\varphi}.
\end{align*}
Then using Young's inequality, we have
\begin{align*}
I_1 & = \sum_j \frac{(1 + z)^\frac{1}{2}}{z} f_{p_j}z_j \notag\\
&\geq - (1 + z)^{\frac{1}{2}} |\nabla_p f| \frac{|\nabla\varphi|}{\varphi} \notag\\
&\geq -\varepsilon (1 + z) |\nabla_p f|^2 - C_{\varepsilon}\frac{|\nabla\varphi|^2}{\varphi^2},
\end{align*}
where $\varepsilon $ is a small constant to be determined.
Note that (\ref{varphi derivative}) implies
\begin{align}\label{grad varphi}
\frac{|\nabla \varphi|^2}{\varphi^2} \leq \frac{C_\alpha}{R^2\varphi^{\frac{2}{\alpha}}}
\leq C_{n,\alpha} I_9.
\end{align}
Then we have
\begin{align*}
I_1 \geq -\varepsilon (1 + z) |\nabla_p f|^2 - C_{n,\alpha,\varepsilon}I_9.
\end{align*}
Using Young's inequality again, $I_7$ can be absorbed by $I_2$ and $I_9$
\begin{align*}
I_7
\leq \frac{z}{n(1 + z)} f^2 + C_{n} \frac{\left| \nabla \varphi \right|^2}{\varphi^2}
\leq \frac{1}{2} I_2 + C_{n,\alpha} I_9.
\end{align*}
Since $F=z$, from (\ref{def G H}) we have
\begin{align}\label{z def G H}
\begin{cases}
\displaystyle
G = -\frac{2+z}{2 z^2(1 + z)^2}<0; \\
\displaystyle
H = \frac{1}{z^2(1 + z)} >0,
\end{cases}
\end{align}
and using (\ref{grad varphi}) again
we obtain the following estimate
\begin{align*}
I_8 \leq -\frac{F^2}{F'^2\varphi^2} G |\nabla\varphi|^2 &
= \frac{(2 + z)|\nabla\varphi|^2}{2(1 + z)^2\varphi^2}
\leq C_{n,\alpha} I_9.
\end{align*}
Substituting all these into (\ref{6e8b74b95}), we arrive at
\begin{align}\label{928e2bde3}
\frac{1 + z}{nz}f^2 -\varepsilon (1 + z) |\nabla_p f|^2
\leq \frac{C_{n,\alpha,\varepsilon}}{R^2\varphi^\frac{2}{\alpha}}.
\end{align}
It follows from our assumptions on $f$ that
\begin{align}\label{bc1c4a907}
\begin{cases}
\theta>0;\\
f^2 - m_1 z |\nabla_p f|^2 \geq m_2 z^\theta.
\end{cases}
\end{align}
We can choose $\varepsilon $ sufficiently small such that
\begin{align*}
\varepsilon \leq \frac{m_1}{n}.
\end{align*}
From (\ref{bc1c4a907}) and (\ref{928e2bde3}), we have
\begin{align}\label{04b4f7260}
z^\theta
\leq C_{m_2} \frac{1 + z}{z} \left( f^2 - m_1 z |\nabla_p f|^2 \right)
\leq C_{n,m_2} \left( \frac{1 + z}{nz} f^2 - \varepsilon (1+z) |\nabla_p f|^2 \right)
\leq \frac{C_{n,m_1,m_2,\alpha}}{R^2\varphi^\frac{2}{\alpha}}.
\end{align}
Let us select $\alpha$ large enough such that
\begin{align*}
\alpha \geq \max\left\{\frac{2}{\theta}, 1 \right\},
\end{align*}
Noting $0< \varphi\leq 1$, then from (\ref{04b4f7260}) we get
\begin{align*}
z\varphi \leq \frac{C_{n,m_1,m_2,\theta}}{R^\frac{2}{\theta}}
\end{align*}
at point $x_0$.
Recalling that $h=1$ and $F=z$.
Since $hF\varphi=z\varphi$ attains its maximum at $x_0$, we have
\begin{align*}
z(0) = (hF\varphi)(0) \leq (hF\varphi)(x_0) \leq \frac{C_{n,m_1,m_2,\theta}}{R^{\frac{2}{\theta}}},
\end{align*}
Using $z = |\nabla u|^2$, we deduce
\begin{align*}
|\nabla u(0)| \leq \frac{C_{n,m_1,m_2,\theta}}{R^{\frac{1}{\theta}}},
\end{align*}
which completes the proof of Theorem \ref{main theorem} (a).

\subsection*{Step 3.}
In this step, we continue from Step 1 to prove Theorem \ref{main theorem} (b).
We set $h=1$, so the following still holds
\begin{align}\label{cf72ad44d}
I_1 + I_2 \leq I_7 + I_8 + I_9.
\end{align}
However, we now set
\begin{align*}
F = \log(1 + z).
\end{align*}
So (\ref{s1 est nabla z}) implies
\begin{align}\label{s3 est nabla z}
\nabla z = -(1 + z)\log(1 + z) \frac{\nabla \varphi}{\varphi}.
\end{align}
Let $\varepsilon $ be a small enough constant to be determined,
using (\ref{grad varphi}) and (\ref{s3 est nabla z}) and Young's inequality,
we have
\begin{align*}
I_1
&\geq -\varepsilon (1 + z) |\nabla_p f|^2 - C_\varepsilon \left| \frac{\nabla\varphi}{\varphi} \right| ^2 \notag\\
&\geq -\varepsilon (1 + z) |\nabla_p f|^2 - \frac{C_{\alpha,\varepsilon}}{R^2\varphi^{\frac{2}{\alpha}}}
\end{align*}
and
\begin{align*}
I_2 = \frac{2f^2}{n\log(1+z)}.
\end{align*}
For $I_7$, we have
\begin{align*}
I_7 \leq \sqrt{\frac{z}{1+z}} f \frac{\left|\nabla\varphi \right|}{\varphi}
\leq C_{\alpha} \sqrt{\frac{z}{1+z}} \frac{f}{R\varphi^{\frac{1}{\alpha}}}.
\end{align*}
Since $F=\log(1+z)$, direct computation gives
\begin{align}\label{logz def G H}
\begin{cases}
\displaystyle
G = \frac{\log(1 + z)-2(1 + z)}{2(1 + z)^3 \log ^2(1 + z)}<0; \\
\displaystyle
H = \frac{1}{(1 + z)^2 \log ^2(1 + z)} > 0.
\end{cases}
\end{align}
Then we have
\begin{align}\label{step3 I8}
I_8 \leq -\frac{F^2}{F'^2\varphi^2} G|\nabla\varphi|^2
\leq \frac{|\nabla\varphi|^2}{\varphi^2}
\leq \frac{C_{\alpha}}{R^2\varphi^{\frac{2}{\alpha}}}.
\end{align}
Hence (\ref{cf72ad44d}) becomes
\begin{align}\label{3ba2544ba}
\frac{2}{n\log(1 + z)} \left( f^2 - \frac{\varepsilon n}{2}(1 + z) \log(1 + z) |\nabla_p f|^2 \right)
&\leq
C_{n,\alpha,\varepsilon}\left( \sqrt{\frac{z}{1+z}} \frac{f}{R\varphi^{\frac{1}{\alpha}}}
+ \frac{1}{R^2\varphi^{\frac{2}{\alpha}}} \right).
\end{align}
Recall that our assumptions on $f$ reads
\begin{align*}
\begin{cases}
\theta > 1;\notag\\
f^2 - m_1 (1+z)\log(1 + z)|\nabla_p f|^2 \geq m_2 \log^{2\theta}(1 + z); \notag\\
|f|<m_3\log^\theta(1 + z).
\end{cases}
\end{align*}
We choose $\varepsilon $ small enough such that
\begin{align*}
\varepsilon < \frac{2m_1}{n},
\end{align*}
then from (\ref{3ba2544ba}) we have
\begin{align*}
\log^{2\theta-1}(1 + z)
& \leq C_{n,m_1,m_2,m_3,\alpha} \left( \log^\theta(1 + z) \frac{1}{R\varphi^{\frac{1}{\alpha}}}
+ \frac{1}{R^2\varphi^{\frac{2}{\alpha}}} \right) \notag\\
&\leq C_{n,m_1,m_2,m_3,\alpha} \max\left\{\log^\theta(1 + z) \frac{1}{R\varphi^{\frac{1}{\alpha}}}, \,
\frac{1}{R^2\varphi^{\frac{2}{\alpha}}} \right\}.
\end{align*}
Let $\alpha$ be large enough such that
\begin{align*}
\alpha \geq \max\left\{\frac{1}{\theta-1},\ 1 \right\}.
\end{align*}
Due to $\theta>1$, observing that
\begin{align}\label{1f87b031a}
\frac{1}{\theta-1} > \frac{2}{2\theta-1},
\end{align}
then we get
\begin{align*}
\log(1 + z)\varphi \leq C_{n,m_1,m_2,m_3,\theta} \max\left\{\frac{1}{R^{\frac{1}{\theta-1}}}, \ \frac{1}{R^{\frac{2}{2\theta-1}}} \right\}.
\end{align*}
Recalling $h=1$ and $F=\log(1+z)$.
Since $hF\varphi = \log(1+z)\varphi$ attains its maximum at $x_0$,
following an argument similar to the proof in Step 2, we have
\begin{align*}
\log(1+z(0)) &= (hF\varphi) (0) \notag\\
&\leq (hF\varphi) (x_0) \notag\\
&\leq C_{n,m_1,m_2,m_3,\theta} \max\left\{\frac{1}{R^{\frac{1}{\theta-1}}}, \ \frac{1}{R^{\frac{2}{2\theta-1}}} \right\}.
\end{align*}
Thus
\begin{align}\label{logz du est 1}
z(0) \leq \max\left\{\exp\left( \frac{C_{n,m_1,m_2,m_3,\theta}}{R^{\frac{1}{\theta-1}}} \right)-1, \ \exp\left( \frac{C_{n,m_1,m_2,m_3,\theta}}{R^{\frac{2}{2\theta-1}}} \right)-1 \right\}.
\end{align}
Using (\ref{1f87b031a}) again, we obtain that, when $R$ is large, the term $\exp\left( CR^{-2/(2\theta-1)} \right)-1$ dominates the right-hand side of (\ref{logz du est 1}).
Thus (\ref{logz du est 1}) can be rewritten as
\begin{align*}
|\nabla u(0)|^2 \leq \max\left\{\exp\left( \frac{C}{R^{\frac{1}{\theta-1}}} \right)-1, \ \frac{C}{R^{\frac{2}{2\theta-1}}} \right\}.
\end{align*}
Here $C$ is a positive constant depending only on $n$, $\theta$, $m_1$, $m_2$ and $m_3$.
This completes the proof.

\subsection*{Step 4.}
In this step, we continue from Step 1 to prove Theorem \ref{main theorem} (c).
Let $F = z$, but define $h$ as follows
\begin{align*}
h = (u + M - 2m + 1)^{-1/b},
\end{align*}
where
\begin{align}\label{def M m}
\begin{cases}
M := \mathop{\sup}\limits_{x\in B_R(0)} u(x), \\
m := \mathop{\inf}\limits_{x\in B_R(0)} u(x).
\end{cases}
\end{align}
and $b$ is a positive constant to be determined later.
The assumption on $f$ reads
\begin{align*}
\begin{cases}
0<\theta\leq 1;\\
f = m_1 \log^{\theta}(1 + z).
\end{cases}
\end{align*}
So there is
\begin{align*}
f_{p_j} = C_{m_1,\theta} \log^{\theta-1}(1 + z) \frac{u_j}{1 + z}.
\end{align*}
Thus using (\ref{s1 est nabla z}) and Young's inequality, we get
\begin{align*}
I_1
& = - C_{m_1,\theta} \frac{\log^{\theta-1}(1 + z)}{(1 + z)^{\frac{1}{2}}} \left( \frac{\nabla u \nabla\varphi}{\varphi} + \frac{h'}{h} z \right) \notag\\
&\geq -\varepsilon _1 \frac{z}{1 + z} \log^{2\theta-2}(1 + z)
- C_{m_1,\theta,\varepsilon _1} \frac{\left| \nabla\varphi \right|^2}{\varphi^2}
- C_{m_1,\theta} \frac{h'}{h} \frac{z}{(1 + z)^{\frac{1}{2}}} \log^{\theta-1}(1 + z),
\end{align*}
where $\varepsilon _1$ is a positive constant to be determined.
and
\begin{align*}
I_2 = \frac{2 m_1^2(1 + z)}{nz}\log^{2\theta}(1 + z).
\end{align*}
It is easy to verify that
\begin{align*}
\frac{1 + z}{z} \log^{2\theta}(1 + z) \geq \frac{z}{1 + z} \log^{2\theta-2}(1 + z),
\qquad \forall z>0.
\end{align*}
This implies
\begin{align*}
I_1 + I_2 \geq &\left( \frac{2m_1^2}{n} - \varepsilon _1 \right) \frac{1 + z}{z} \log^{2\theta}(1 + z) \notag\\
&- C_{m_1,\theta,\varepsilon _1} \frac{\left| \nabla\varphi \right|^2}{\varphi^2}
- C_{m_1,\theta} \frac{h'}{h} \frac{z}{(1 + z)^{\frac{1}{2}}} \log^{\theta -1}(1 + z).
\end{align*}
Using Young's inequality, we have
\begin{align*}
I_3 & = m_1 \frac{h'}{h} \frac{\log^{\theta}(1 + z)}{(1 + z)^{\frac{1}{2}}} \notag\\
&\geq - \frac{\varepsilon _2}{z} \log^{2\theta}(1 + z)
- \frac{m_1^2}{4\varepsilon _2} \left( \frac{h'}{h} \right)^2 \frac{z}{1 + z} \notag\\
&\geq - \varepsilon _2\frac{1 + z}{z} \log^{2\theta}(1 + z)
- \frac{m_1^2}{4\varepsilon _2} \left( \frac{h'}{h} \right)^2 \frac{z}{1 + z},
\end{align*}
where $\varepsilon _2$ is a positive constant to be determined.
Note that $h$ satisfies
\begin{align}\label{h prop}
\frac{h''}{h} = (b + 1)\frac{h'^2}{h^2}.
\end{align}
Thus we have
\begin{align*}
I_4 = b \left( \frac{h'}{h} \right)^2 \frac{z}{1 + z}.
\end{align*}
Using the same values of $G$ and $H$ in (\ref{z def G H}), we obtain
\begin{align*}
I_5 = \frac{1}{2} \left( \frac{h'}{h} \right)^2 \frac{z(z-2)}{(1 + z)^2},
\end{align*}
and
\begin{align*}
I_6 & = -\frac{h'}{h} \frac{z-2}{(1 + z)^2} \frac{\nabla u \nabla \varphi}{\varphi} \notag\\
& \leq \frac{1}{8} \left( \frac{h'}{h} \right)^2 \frac{z(z-2)^2}{(1 + z)^4}
+ C \frac{|\nabla \varphi|^2}{\varphi^2}.
\end{align*}
Noting for all $z\geq 0$ there is
\begin{align*}
-2 \leq \frac{z-2}{(1 + z)^{2}} \leq \frac{1}{12}.
\end{align*}
This lead us to set
\begin{align*}
\kappa =
\begin{cases}
-2, & \text{if}\ 0 \leq z \leq 2; \notag\\
\frac{1}{12}, & \text{if}\ z>2.
\end{cases}
\end{align*}
Then we have
\begin{align*}
I_6 \leq \frac{\kappa}{8} \left( \frac{h'}{h} \right)^2 \frac{z(z-2)}{(1 + z)^2}
+ C \frac{|\nabla \varphi|^2}{\varphi^2}.
\end{align*}
Furthermore
\begin{align*}
I_5 - I_6 &\geq \frac{4-\kappa}{8} \left( \frac{h'}{h} \right)^2 \frac{z(z-2)}{(1 + z)^2}
- C \frac{|\nabla \varphi|^2}{\varphi^2} \notag\\
&\geq -\frac{4-\kappa}{4} \left( \frac{h'}{h} \right)^2 \frac{z}{1 + z}
- C \frac{|\nabla \varphi|^2}{\varphi^2}.
\end{align*}
Direct computation gives
\begin{align*}
I_7 & = m_1 \frac{\log^{\theta}(1 + z)}{(1 + z)^{\frac{1}{2}}} \frac{\nabla u \nabla \varphi}{\varphi} \notag\\
& \leq \varepsilon _3\frac{z}{(1 + z)} \log^{2\theta}(1 + z)
+ C_{m_1,\varepsilon _3} \frac{\left| \nabla\varphi \right|^2}{\varphi^2} \notag\\
& \leq \varepsilon _3\frac{1 + z}{z} \log^{2\theta}(1 + z)
+ C_{m_1,\varepsilon _3} \frac{\left| \nabla\varphi \right|^2}{\varphi^2},
\end{align*}
where $\varepsilon _3$ is a positive constant to be determined.
Using (\ref{z def G H}) again,
$I_8$ is estimated as follows
\begin{align*}
I_8 \leq - \frac{F^2G|\nabla\varphi|^2}{F'^2\varphi^2}
\leq \frac{|\nabla \varphi|^2}{\varphi^2}.
\end{align*}
Combining these estimates, and using (\ref{grad varphi}), then (\ref{step1 result}) becomes
\begin{align}\label{2ab2496a2}
\left( \frac{2m_1^2}{n} - \varepsilon _1 -\varepsilon _2 -\varepsilon _3 \right) \frac{1 + z}{z} \log^{2\theta}(1 + z) &\notag\\
+ \left( b-\frac{m_1^2}{4\varepsilon _2} - \frac{4-\kappa}{4} \right) \left( \frac{h'}{h} \right)^2 \frac{z}{1 + z}& \notag\\
- C_{m_1,\theta} \frac{h'}{h} \frac{z}{(1 + z)^{\frac{1}{2}}} \log^{\theta-1}(1 + z)
&\leq \frac{C_{m_1,\alpha,\theta,\varepsilon _1,\varepsilon _3}}{R^2\varphi^{\frac{2}{\alpha}}}.
\end{align}
We choose $\varepsilon _1$, $\varepsilon _2$ and $\varepsilon _3$ small enough such that
\begin{align*}
\varepsilon _1 + \varepsilon _2 + \varepsilon _3 \leq \frac{2m_1^2}{n}.
\end{align*}
Next, we choose $b$ large enough such that
\begin{align*}
b \geq \frac{m_1^2}{4\varepsilon _2} + \frac{4-\kappa}{4}.
\end{align*}
Noting that $h' = -\frac{1}{b} h^{b+1}$, then (\ref{2ab2496a2}) implies that
\begin{align}\label{7c01d5014}
\frac{h^bz}{(1 + z)^{\frac{1}{2}}} \log^{\theta-1}(1 + z)
\leq
\frac{C_{n,m_1,\alpha,\theta}}{R^2\varphi^{\frac{2}{\alpha}}}.
\end{align}

We now divide the proof into two cases:

\textbf{case 1:}
For all small $z$, there exists constant $C$ such that $1 + z \leq C$.
Using the elementary inequality $\log(1+z)\leq z$,
form (\ref{7c01d5014}), we obtain
\begin{align}\label{f09293cf3}
h^b z^{\theta}
\leq \frac{Ch^bz}{(1 + z)^{\frac{1}{2}}} \log^{\theta-1}(1 + z)
\leq \frac{C_{n,m_1,\alpha,\theta}}{R^2\varphi^{\frac{2}{\alpha}}}
\end{align}
and we may set $\alpha$ large enough such that
\begin{align*}
\alpha \geq \frac{2}{\theta},
\end{align*}
then (\ref{f09293cf3}) implies
\begin{align*}
hF\varphi \leq
C_{n,m_1,\theta}\frac{h^{1-\frac{b}{\theta}}}{R^{\frac{2}{\theta}}},
\end{align*}
which complete the case 1.

\textbf{case 2:}
For all large $z$, there exists constant $C_\eta$ and $C$ such that
\begin{align*}
\begin{cases}
\log^{1-\theta}(1 + z) \leq C_\eta z^{\frac{\eta}{2}}; \\
1 + z \leq Cz.
\end{cases}
\end{align*}
Then (\ref{7c01d5014}) implies
\begin{align*}
h^b z^{\frac{1-\eta}{2}} \leq \frac{C_{n,m_1,\alpha,\eta,\theta}}{R^2\varphi^{\frac{2}{\alpha}}}.
\end{align*}
We choose $\alpha$ large enough such that
\begin{align*}
\alpha \geq \frac{4}{1-\eta},
\end{align*}
then we obtain
\begin{align*}
hF\varphi \leq C_{n,m_1,\eta,\theta} \frac{h^{1-\frac{2b}{1-\eta}}}{R^\frac{4}{1-\eta}},
\end{align*}
which complete the case 2.

We may choose $\alpha$ large enough such that
\begin{align*}
\alpha \geq \max\left\{\frac{2}{\theta}, \ \frac{4}{1-\eta} \right\}.
\end{align*}
Then combining these two cases, we have
\begin{align*}
hF\varphi \leq \max \left\{
C_{n,m_1,\theta} \frac{h^{1-\frac{b}{\theta}}}{R^{\frac{2}{\theta}}}, \
C_{n,m_1,\eta,\theta} \frac{h^{1-\frac{2b}{1-\eta}}}{R^\frac{4}{1-\eta}}
\right\}
\end{align*}
at point $x_0$.
Similar to the proof in the previous step, we have
\begin{align*}
h(0)z(0) &=
(hF\varphi) (0) \notag\\
&\leq (hF\varphi) (x_0)\notag\\
&\leq \max \left\{
C_{n,m_1,\theta}\frac{h^{1-\frac{b}{\theta}}(x_0)}{R^{\frac{2}{\theta}}}, \
C_{n,m_1,\eta,\theta}\frac{h^{1-\frac{2b}{1-\eta}}(x_0)}{R^\frac{4}{1-\eta}}
\right\}.
\end{align*}
Note that $h$ satisfies the following two-sided estimate:
\begin{align*}
\frac{1}{2^{\frac{1}{b}}(L + 1)^{\frac{1}{b}}} \leq h(x) \leq \frac{1}{(L + 1)^{\frac{1}{b}}}, \quad \forall x\in B_R(0).
\end{align*}
We will finally arrive at
\begin{align*}
|\nabla u(0)| \leq \max\left\{
C_1\frac{(L + 1)^\frac{1}{2\theta}}{R^\frac{1}{\theta}}, \
C_2\frac{(L + 1)^\frac{1}{1 - \eta}}{R^\frac{2}{1-\eta}}
\right\}.
\end{align*}
Here $C_1$ is a positive constant depending on $n$, $\theta$ and $m_1$.
$C_2$ is a positive constant depending on $n$, $\eta$, $\theta$ and $m_1$.
This complete the proof of Theorem \ref{main theorem} (c).

\subsection*{Step 5.}
In this step, we continue from Step 1 to prove theorem \ref{main theorem} (d).
Let
\begin{align}\label{2b9f7b8b2}
\begin{cases}
F = \log(1 + z); \\
h = (u + M-2m+1)^{-1/b}.
\end{cases}
\end{align}
Here, $M$ and $m$ are defined the same as in (\ref{def M m}), and $b$ is a constant to be determined.
Using (\ref{s1 est nabla z}), we get
\begin{align*}
I_1 & \geq -(1 + z)^{\frac{1}{2}} |\nabla_p f|\left( \frac{|\nabla \varphi|}{\varphi} + \frac{|h'|}{h}|\nabla u| \right) \notag\\
&\geq -\varepsilon _1(1 + z)^2 |\nabla_p f|^2
- \frac{1}{4\varepsilon _1} \left( \frac{h'}{h} \right)^2 \frac{z}{1 + z}
- \frac{C_{\varepsilon _1}}{1 + z} \frac{\left| \nabla\varphi \right|^2}{\varphi^2},
\end{align*}
where $\varepsilon _1$ is positive constant to be determined.
Direct computation shows
\begin{align*}
I_2 = \frac{2f^2}{n\log(1+z)}.
\end{align*}
Using Young's inequality and inequality $\log(1+z)\leq z$, we have
\begin{align*}
I_3 &\geq - \frac{1}{nz} f^2 - \frac{n}{4} \left( \frac{h'}{h} \right)^2 \frac{z}{1 + z}\notag\\
&\geq -\frac{1}{n\log(1 + z)} f^2 - \frac{n}{4} \left( \frac{h'}{h} \right)^2 \frac{z}{1 + z}.
\end{align*}
The assumptions on $f$ read as follows:
\begin{align}\label{step5 assumption f}
\begin{cases}
f^2 - m_1(1 + z)^2\log(1 + z) |\nabla_p f|^2 \geq 0; \\
|f| \leq m_2.
\end{cases}
\end{align}
We deduce that
\begin{align*}
I_1 + I_2 + I_3
\geq & \frac{1}{n\log(1+z)} \left( f^2 - n\varepsilon _1 (1+z)^2\log(1+z) |\nabla_p f|^2 \right) \notag\\
& -\left( \frac{1}{4\varepsilon _1} + \frac{n}{4} \right) \left( \frac{h'}{h} \right)^2 \frac{z}{1 + z}
- \frac{C_{n,m_1}}{1 + z} \frac{\left| \nabla\varphi \right|^2}{\varphi^2}.
\end{align*}
We choose $\varepsilon _1$ small enough such that
\begin{align*}
\varepsilon _1 \leq \frac{m_1}{n},
\end{align*}
then we get
\begin{align*}
I_1 + I_2 + I_3
\geq
-\left( \frac{1}{4\varepsilon _1} + \frac{n}{4} \right) \left( \frac{h'}{h} \right)^2 \frac{z}{1 + z}
- \frac{C_{n,m_1}}{1 + z} \frac{\left| \nabla\varphi \right|^2}{\varphi^2}.
\end{align*}
Noting that in this step $h$ also satisfies (\ref{h prop}), we have
\begin{align*}
I_4 = b \left( \frac{h'}{h} \right)^2 \frac{z}{1 + z}.
\end{align*}
Direct computation gives
\begin{align}\label{logz main term est}
I_5
= \frac{1}{2} \left( \frac{h'}{h} \right)^2 \frac{z}{1 + z}\Big( \log(1 + z)-2 \Big).
\end{align}
Using the same value of $G$ and $H$ in (\ref{logz def G H}) and Young's inequality, we get
\begin{align*}
I_6 & \leq \frac{|h'|}{h} \frac{z^{\frac{1}{2}}(\log(1 + z)-2)}{1 + z} \frac{|\nabla \varphi|}{\varphi} \notag\\
& \leq \frac{1}{8} \left( \frac{h'}{h} \right)^2 \frac{z}{(1 + z)^2} \Big( \log(1 + z)-2 \Big)^2
+ C\frac{|\nabla \varphi|^2}{\varphi^2}.
\end{align*}
Since
\begin{align}
-2 \leq \frac{\log(1 + z)-2}{1 + z} \leq e^{-3}, \quad \forall z\geq 0,
\end{align}
we may set
\begin{align*}
\kappa = \left\{
\begin{array}{ll}
e^{-3}, \quad& \text{if}\ \log(1 + z)>2; \\
-2, & \text{if}\ \log(1 + z)<2.
\end{array}
\right.
\end{align*}
Then we get
\begin{align*}
I_6 \leq \frac{\kappa}{8} \left( \frac{h'}{h} \right)^2 \frac{z}{1 + z}\Big( \log(1 + z)-2 \Big)
+ C \frac{|\nabla \varphi|^2}{\varphi^2}.
\end{align*}
Combining with (\ref{logz main term est}), we get
\begin{align*}
I_5-I_6
\geq \frac{4-\kappa}{8} \left( \frac{h'}{h} \right)^2 \frac{z}{1 + z}\Big( \log(1 + z)-2 \Big)
- C\frac{|\nabla \varphi|^2}{\varphi^2}.
\end{align*}
Then we have
\begin{align*}
I_1 + I_2 + I_3 + I_4 + I_5 - I_6
\geq &\left( b - \frac{1}{4\varepsilon _1} - \frac{n}{4} - \frac{4-\kappa}{4} \right)\left( \frac{h'}{h} \right)^2 \frac{z}{1 + z} \notag\\
&+ \frac{4-\kappa}{8} \left( \frac{h'}{h} \right)^2 \frac{z}{1 + z} \log(1 + z)
- C_{n,m_1} \frac{|\nabla \varphi|^2}{\varphi^2}.
\end{align*}
We choose $b$ large enough such that
\begin{align*}
b\geq \frac{1}{4\varepsilon _1} + \frac{n}{4} + \frac{4-\kappa}{4}.
\end{align*}
Thus we have
\begin{align*}
I_1 + I_2 + I_3 + I_4 + I_5 - I_6 \geq \frac{1}{4} \left( \frac{h'}{h} \right)^2 \frac{z}{1 + z} \log(1 + z)
- C_{n,m_1} \frac{|\nabla \varphi|^2}{\varphi^2}.
\end{align*}
Direct computation gives
\begin{align*}
I_7 &\leq m_2 \sqrt{\frac{z}{1 + z}} \left| \frac{\nabla\varphi}{\varphi} \right|.
\end{align*}
Utilizing (\ref{logz def G H}) again,
we can obtain the same estimate for $I_8$ as in (\ref{step3 I8}), that is
\begin{align*}
I_8 \leq \frac{|\nabla\varphi|^2}{\varphi^2}.
\end{align*}
Then by (\ref{grad varphi}), from (\ref{step1 result}) we obtain
\begin{align}\label{logz est 1}
\left( \frac{h'}{h} \right)^2 \frac{z}{1 + z} \log(1 + z)
\leq C_{n,m_1,m_2,\alpha} \left( \sqrt{\frac{z}{1 + z}} \frac{1}{R\varphi^\frac{1}{\alpha}}
+ \frac{1}{R^2\varphi^\frac{2}{\alpha}}
\right).
\end{align}

Let $\beta$ and $\gamma$ be undetermined constants satisfying $0<\beta\leq 2$ and $\gamma>0$.
We split the proof into two cases:

\textbf{cases 1:} Assume that $z$ is sufficiently small at point $x_0$ such that
\begin{align*}
z \leq \frac{h^{-\gamma}}{\varphi^{\frac{\beta}{\alpha}}R^{\beta}}.
\end{align*}
Then by the inequality $\log(1 + z) \leq z$, we get
\begin{align*}
\log(1 + z) \leq \frac{h^{-\gamma}}{\varphi^{\frac{\beta}{\alpha}}R^\beta}.
\end{align*}
We choose $\alpha$ large enough such that
\begin{align}\label{8bcb2934c}
\alpha \geq 2.
\end{align}
Recalling $\varphi \leq 1$ and noting
\begin{align}\label{two-sided est h}
\frac{1}{2^{\frac{1}{b}}(L + 1)^{\frac{1}{b}}} \leq h \leq \frac{1}{(L + 1)^{\frac{1}{b}}},
\end{align}
we get
\begin{align}\label{case1 main ret}
\log(1 + z)\varphi \leq C_{n,m_1}\frac{(L + 1)^{\frac{\gamma}{b}}}{R^\beta},
\end{align}
which complete the case 1.

\textbf{cases 2:} Otherwise, we have
\begin{align*}
z > \frac{h^{-\gamma}}{\varphi^{\frac{\beta}{\alpha}}R^{\beta}}
\end{align*}
at point $x_0$.
Then use the elementary inequality $\sqrt{1 + t} \leq 1 + \sqrt{t} $, from (\ref{logz est 1}) we have
\begin{align}\label{290e50c8f}
\left( \frac{h'}{h} \right)^2 \log(1 + z)
& \leq C_{n,m_1,m_2,\alpha}\left( ( 1 + \frac{1}{\sqrt{z}} ) \frac{1}{R\varphi^\frac{1}{\alpha}}
+ (1 + \frac{1}{z} )\frac{1}{R^2\varphi^\frac{2}{\alpha}} \right) \notag\\
& \leq C_{n,m_1,m_2,\alpha}\left( \frac{1}{R \varphi^{\frac{1}{\alpha}}}
+ \frac{h^{\frac{\gamma}{2}}}{R^{1-\frac{\beta}{2}} \varphi^{\frac{2-\beta}{2\alpha}}}
+ \frac{1}{R^2 \varphi^{\frac{2}{\alpha}}}
+ \frac{h^{\gamma}}{R^{2-\beta} \varphi^{\frac{2-\beta}{\alpha}}}
\right).
\end{align}
In this case we choose the same $\alpha$ satisfying (\ref{8bcb2934c}) as in case 1.
Noting $h' =-\frac{1}{b} h^{b+1}$,
then (\ref{290e50c8f}) gives
\begin{align*}
\log(1 + z)\varphi \leq C_{n,m_1,m_2} \max\left\{
\frac{h^{-2b}}{R}, \
\frac{h^{\frac{\gamma}{2}-2b}}{R^{1-\frac{\beta}{2}}}, \
\frac{h^{-2b}}{R^2}, \
\frac{h^{\gamma-2b}}{R^{2-\beta}}
\right\}.
\end{align*}
Then using (\ref{two-sided est h}) again, we have
\begin{align}\label{case2 main ret}
\log(1 + z)\varphi & \leq C_{n,m_1,m_2} \max\left\{
\frac{(L + 1)^2}{R}, \
\frac{(L + 1)^{2-\frac{\gamma}{2b}}}{R^{1-\frac{\beta}{2}}}, \
\frac{(L + 1)^2}{R^2}, \
\frac{(L + 1)^{2-\frac{\gamma}{b}}}{R^{2-\beta}}
\right\},
\end{align}
which complete the case 2.

Combining the two cases (\ref{case1 main ret}) and (\ref{case2 main ret}), we have
\begin{align}\label{s5 combining two cases}
&\log(1 + z)\varphi \notag\\
&\leq C_{n,m_1,m_2} \max\left\{
\frac{(L + 1)^{\frac{\gamma}{b}}}{R^\beta}, \
\frac{(L + 1)^2}{R}, \
\frac{(L + 1)^{2-\frac{\gamma}{2b}}}{R^{1-\frac{\beta}{2}}}, \
\frac{(L + 1)^2}{R^2}, \
\frac{(L + 1)^{2-\frac{\gamma}{b}}}{R^{2-\beta}}
\right\}.
\end{align}
Let $\gamma$ small enough such that
\begin{align*}
\gamma < 2b,
\end{align*}
then we have
\begin{align}\label{logz log leq 5 terms}
\log(1 + z)\varphi \leq C_{n,m_1,m_2} (L + 1)^{2} \max\left\{
\frac{1}{R^\beta}, \
\frac{1}{R}, \
\frac{1}{R^{1-\frac{\beta}{2}}}, \
\frac{1}{R^2}, \
\frac{1}{R^{2-\beta}}
\right\}.
\end{align}
Define the set $\mathbb{S}$ as the collection of exponents of $R$ in the right-hand side of above inequality, i.e.,
\begin{align*}
\mathbb{S} &:= \left\{
\beta, \ 1, \ {1-\frac{\beta}{2}}, \ 2, \ 2-\beta
\right\} \notag\\
\end{align*}
and define
\begin{align*}
\lambda_1 &:= \inf\limits_{0<\beta\leq 2} \max \mathbb{S}, \notag\\
\lambda_2 &:= \sup\limits_{0<\beta\leq 2} \min \mathbb{S}.
\end{align*}

When $R<1$, we can choose a suitable $\beta$ such that the right-hand side of (\ref{logz log leq 5 terms}) is controlled by $(L + 1)^{2}/R^{\lambda_1}$;
on the other hand, when $R\geq 1$, it will be controlled by $(L + 1)^{2}/R^{\lambda_2}$. Combining these two cases, we obtain
\begin{align*}
\log(1 + z)\varphi & \leq
C_{n,m_1,m_2} (L + 1)^{2} \max\left\{\frac{1}{R^{\lambda_1}}, \ \frac{1}{R^{\lambda_2}} \right\}.
\end{align*}
Recalling (\ref{2b9f7b8b2}), we conclude
\begin{align*}
h(0) \log(1+z(0)) &= (hF\varphi) (0) \notag\\
&\leq (hF\varphi) (x_0) \notag\\
&\leq C_{n,m_1,m_2} h(x_0) \max\left\{\frac{(L + 1)^{2}}{R^{\lambda_1}}, \ \frac{(L + 1)^{2}}{R^{\lambda_2}} \right\}.
\end{align*}
Using (\ref{two-sided est h}) again, we get the gradient estimate at the center of the ball $B_R(0)$ is of the form
\begin{align*}
|\nabla u(0)|^2 \leq \max\left\{\exp\left( C_{n,m_1,m_2} \frac{(L + 1)^2}{R^{\lambda_1}} \right) - 1, \
\exp\left( C_{n,m_1,m_2} \frac{(L + 1)^{2}}{R^{\lambda_2}}\right)-1
\right\}.
\end{align*}
Finally, we select an appropriate $\beta$ to determine $\lambda_1$ and $\lambda_2$.
Direct computation shows that
the infimum of $\max \mathbb{S}$ is realized when $0\leq \beta \leq 2$, yielding $\lambda_1=2$;
the supremun of $\min \mathbb{S}$ is realized when $\beta = 2/3$, yielding $\lambda_2=2/3$.
Thus we get
\begin{align}\label{4725fc0a5}
|\nabla u(0)|^2 \leq \max
\left\{\exp\left( C_{n,m_1,m_2}\frac{(L + 1)^2}{R^2} \right) - 1,
\exp\left( C_{n,m_1,m_2}\frac{(L + 1)^{2}}{R^{\frac{2}{3}}}\right)-1 \right\}.
\end{align}
which complete the proof of Theorem \ref{main theorem} (d).

\subsection*{Step 6.}
In this step, we prove theorem \ref{main theorem} (e).
Noting that $f=0$, which satisfies (\ref{step5 assumption f}),
thus we may adopt nearly the same proof as in Step 5.
However, in contrast, we set
\begin{align*}
h=(u+M-2m)^{-1/b}.
\end{align*}
Thus (\ref{two-sided est h}) is replaced by
\begin{align*}
\frac{1}{2^{\frac{1}{b}}L^{\frac{1}{b}}} \leq h \leq \frac{1}{L^{\frac{1}{b}}}.
\end{align*}
Moreover $I_7=0$, which implies that (\ref{logz est 1}) reduces to
\begin{align*}
\left( \frac{h'}{h} \right)^2 \frac{z}{1 + z} \log(1 + z)
\leq \frac{C_{n,\alpha}}{R^2\varphi^\frac{2}{\alpha}}.
\end{align*}
After discussing the same two cases, we find that at the point $x_0$, (\ref{s5 combining two cases}) simplifies into
\begin{align*}
\log(1 + z)\varphi & \leq C_{n}\max\left\{
\frac{L^{\frac{\gamma}{b}}}{R^\beta}, \
\frac{L^2}{R^2}, \
\frac{L^{2-\frac{\gamma}{b}}}{R^{2-\beta}}
\right\}.
\end{align*}
Now we can directly set
\begin{align*}
\begin{cases}
\beta=1;\\
\gamma=b.
\end{cases}
\end{align*}
Then following similar proof as in previous step, we derive an enhanced version of the gradient estimate in (\ref{4725fc0a5}), which is
\begin{align*}
|\nabla u(0)|^2 \leq \max\left\{\exp\left( \frac{CL^2}{R^2} \right)-1,\ \exp\left( \frac{CL}{R} \right)-1 \right\}.
\end{align*}
Here $C$ is a positive constant depending only on $n$.
Then we complete the proof of Theorem \ref{main theorem} (e).

\textbf{Acknowledgements.}
The author would like to thank the anonymous referees for their careful reading
of the manuscript and for many valuable comments and suggestions, which have
led to a substantial improvement of the paper.

\textbf{Conflict of Interest Statement.}
The author declares that there is no conflict of interest.

\textbf{Data Availability Statement.}
No datasets were generated or analysed during the current study.


\begin{thebibliography}{99}
\bibitem{Almgren Annals 1966}F.~J. Almgren Jr., Some interior regularity theorems for minimal surfaces and an extension of Bernstein's theorem, Ann. of Math. (2) {\bf 84} (1966), 277--292; MR0200816
\bibitem{Bernstein MathAnn 1906}S. Bernstein, Sur la g\'en\'eralisation du probl\`eme de Dirichlet, Math. Ann. {\bf 62} (1906), no.~2, 253--271; MR1511375
\bibitem{Bernstein MathAnn 1910}S. Bernstein, Sur la g\'en\'eralisation du probl\`eme de Dirichlet, Math. Ann. {\bf 69} (1910), no.~1, 82--136; MR1511579
\bibitem{Bernstein MathZ 1927}S. Bernstein, \"Uber ein geometrisches Theorem und seine Anwendung auf die partiellen Differentialgleichungen vom elliptischen Typus, Math. Z. {\bf 26} (1927), no.~1, 551--558; MR1544873
\bibitem{Bianchini Colombo Magliaro Mari Pucci Rigoli 2021 MathEng}B. Bianchini, G. Colombo, M. Magliaro, L. Mari, P. Pucci, M. Rigoli, Recent rigidity results for graphs with prescribed mean curvature, Math. Eng. {\bf 3} (2021), no.~5, Paper No. 039, 48 pp.; MR4181196
\bibitem{Bianchini Mari Pucci Rigoli 2021 Springer}B. Bianchini, L, Mari, P. Pucci, M. Rigoli, {\it Geometric analysis of quasilinear inequalities on complete manifolds---maximum and compact support principles and detours on manifolds}, Frontiers in Mathematics, Birkh\"auser/Springer, Cham, [2021] \copyright 2021; MR4241012
\bibitem{Bidaut-Veron Gacia-Huidobro Veron Duke 2019}M.-F. Bidaut-V\'eron, M. Garc\'ia-Huidobro and L. V\'eron, Estimates of solutions of elliptic equations with a source reaction term involving the product of the function and its gradient, Duke Math. J. {\bf 168} (2019), no.~8, 1487--1537; MR3959864
\bibitem{Bombieri Giorgi Giusti Invent 1969}E. Bombieri, E. De~Giorgi and E. Giusti, Minimal cones and the Bernstein problem, Invent. Math. {\bf 7} (1969), 243--268; MR0250205
\bibitem{Bombieri Giorgi Miranda ARMA 1969}E. Bombieri, E. De~Giorgi and M. Miranda, Una maggiorazione a priori relativa alle ipersuperfici minimali non parametriche, Arch. Rational Mech. Anal. {\bf 32} (1969), 255--267; MR0248647
\bibitem{Bombieri Giusti Invent 1972}E. Bombieri and E. Giusti, Harnack's inequality for elliptic differential equations on minimal surfaces, Invent. Math. {\bf 15} (1972), 24--46; MR0308945
\bibitem{Cabre Dipierro Valdinoci ARMA 2022}X. Cabr\'e, S. Dipierro and E. Valdinoci, The Bernstein technique for integro-differential equations, Arch. Ration. Mech. Anal. {\bf 243} (2022), no.~3, 1597--1652; MR4381148
\bibitem{Caffarelli Garofalo Segala CPAM 1994}L.~\'A. Caffarelli, N. Garofalo and F. Seg\`ala, A gradient bound for entire solutions of quasi-linear equations and its consequences, Comm. Pure Appl. Math. {\bf 47} (1994), no.~11, 1457--1473; MR1296785
\bibitem{Chou Wang CPAM 2001}K.-S. Chou and X.~J. Wang, A variational theory of the Hessian equation, Comm. Pure Appl. Math. {\bf 54} (2001), no.~9, 1029--1064; MR1835381
\bibitem{Cozzi Farina Valdinoci CMP 2014}M. Cozzi, A. Farina and E. Valdinoci, Gradient bounds and rigidity results for singular, degenerate, anisotropic partial differential equations, Comm. Math. Phys. {\bf 331} (2014), no.~1, 189--214; MR3231999
\bibitem{DeGiorgi ASNSP 1965}E. De~Giorgi, Una estensione del teorema di Bernstein, Ann. Scuola Norm. Sup. Pisa Cl. Sci. (3) {\bf 19} (1965), 79--85; MR0178385
\bibitem{DeSilva Jerison CAPM 2011}D. De~Silva and D.~S. Jerison, A gradient bound for free boundary graphs, Comm. Pure Appl. Math. {\bf 64} (2011), no.~4, 538--555; MR2796515
\bibitem{Ecker Huisken JDG 1990} K. Ecker and G. Huisken, A Bernstein result for minimal graphs of controlled growth, J. Differential Geom. {\bf 31} (1990), no.~2, 397--400; MR1037408
\bibitem{Farina Valdinoci CVPDE 2014}A. Farina and E. Valdinoci, Gradient bounds for anisotropic partial differential equations, Calc. Var. Partial Differential Equations {\bf 49} (2014), no.~3-4, 923--936; MR3168616
\bibitem{Finn ARMA 1963}R.~S. Finn, New estimates for equations of minimal surface type, Arch. Rational Mech. Anal. {\bf 14} (1963), 337--375; MR0157096
\bibitem{Fleming RCMP 1962}W.~H. Fleming, On the oriented Plateau problem, Rend. Circ. Mat. Palermo (2) {\bf 11} (1962), 69--90; MR0157263
\bibitem{Gilbarg Trudinger Berlin 1998}D. Gilbarg and N.~S. Trudinger, {\it Elliptic partial differential equations of second order}, reprint of the 1998 edition, Classics in Mathematics, Springer, Berlin, 2001; MR1814364
\bibitem{Huisken Ilmanen JDG 2001}G. Huisken and T. Ilmanen, The inverse mean curvature flow and the Riemannian Penrose inequality, J. Differential Geom. {\bf 59} (2001), no.~3, 353--437; MR1916951
\bibitem{Jenkins Serrin 1968 JRAM}H. Jenkins and J.~B. Serrin Jr., The Dirichlet problem for the minimal surface equation in higher dimensions, J. Reine Angew. Math. {\bf 229} (1968), 170--187; MR0222467
\bibitem{Korevaar PSPM 1983}N.~J. Korevaar, An easy proof of the interior gradient bound for solutions to the prescribed mean curvature equation, in {\it Nonlinear functional analysis and its applications, Part 2 (Berkeley, Calif., 1983)}, 81--89, Proc. Sympos. Pure Math., 45, Part 2, Amer. Math. Soc., Providence, RI, ; MR0843597
\bibitem{Ladyzhenskaya Uraltseva book 1968}O.~A. Ladyzhenskaya and N.~N. Ural'tseva, {\it Linear and quasilinear elliptic equations}, translated from the Russian by Scripta Technica, Inc Translation editor: Leon Ehrenpreis, Academic Press, New York-London, 1968; MR0244627
\bibitem{Ladyzhenskaya Uraltseva CPAM 1970}O.~A. Ladyzhenskaya and N.~N. Ural'tseva, Local estimates for gradients of solutions of non-uniformly elliptic and parabolic equations, Comm. Pure Appl. Math. {\bf 23} (1970), 677--703; MR0265745
\bibitem{Ma Xu Adv 2016} X.~N. Ma and J. Xu, Gradient estimates of mean curvature equations with Neumann boundary value problems, Adv. Math. {\bf 290} (2016), 1010--1039; MR3451945
\bibitem{Moser CPAM 1961}J.~K. Moser, On Harnack's theorem for elliptic differential equations, Comm. Pure Appl. Math. {\bf 14} (1961), 577--591; MR0159138
\bibitem{Serrin JFA 1970}J.~B. Serrin Jr., On the strong maximum principle for quasilinear second order differential inequalities, J. Functional Analysis {\bf 5} (1970), 184--193; MR0259328
\bibitem{Serrin PTRSLA 1969}J.~B. Serrin Jr., The problem of Dirichlet for quasilinear elliptic differential equations with many independent variables, Philos. Trans. Roy. Soc. London Ser. A {\bf 264} (1969), 413--496; MR0282058
\bibitem{Simon Indiana 1976}L.~M. Simon, Interior gradient bounds for non-uniformly elliptic equations, Indiana Univ. Math. J. {\bf 25} (1976), no.~9, 821--855; MR0412605
\bibitem{Simon JDG 1989} L.~M. Simon, Entire solutions of the minimal surface equation, J. Differential Geom. {\bf 30} (1989), no.~3, 643--688; MR1021370
\bibitem{Simons Annals 1968}J.~H. Simons, Minimal varieties in riemannian manifolds, Ann. of Math. (2) {\bf 88} (1968), 62--105; MR0233295
\bibitem{Trudinger PNAS 1972}N.~S. Trudinger, A new proof of the interior gradient bound for the minimal surface equation in $n$ dimensions, Proc. Nat. Acad. Sci. U.S.A. {\bf 69} (1972), 821--823; MR0296832
\bibitem{Wang MathZ 1998}X.~J. Wang, Interior gradient estimates for mean curvature equations, Math. Z. {\bf 228} (1998), no.~1, 73--81; MR1617971

\end{thebibliography}
\end{document}